\documentclass[12pt]{article}
\usepackage{pgf,tikz}
\usepackage{calc,graphicx}
\usepackage{graphics}
\everymath{\displaystyle}
\usepackage{graphicx}
\usepackage{color}
\usepackage{amssymb}
\usepackage{amsmath}
\newtheorem{prethm}{{\bf Theorem}}

\newenvironment{thm}{\begin{prethm}{\hspace{-0.5
				em}{\bf .}}}{\end{prethm}}
\newtheorem{prelemma}{{\bf Lemma}}

\newtheorem{preex}{{\bf Example}}

\newtheorem{preprop}{{\bf Proposition}}

\newenvironment{prop}{\begin{preprop}{\hspace{-0.5em}{\bf .}}}{\end{preprop}}
\newtheorem{precor}{{\bf Corollary}}

\newenvironment{cor}{\begin{precor}{\hspace{-0.5
				em}{\bf .}}}{\end{precor}}
\newtheorem{preremark}{{\bf Remark}}

\newtheorem{preprob}{{\bf Problem}}

\newtheorem{predefin}{{\bf Definition}}

\newenvironment{defin}{\begin{predefin}{\hspace{-0.5
				em}{\bf .}}}{\end{predefin}}
\newtheorem{preconj}{{\bf Conjecture}}

\newtheorem{preprobb}{{\bf Problem}}

\newtheorem{prelem}{{\bf Theorem}}

\newenvironment{proof}{{\bf Proof.}\rm }{\hfill{$\Box$}}

\newtheorem{presolution}{{\bf Solution.}}

\def\newpic#1{}

\title{\vspace{-0.51cm}\Large\bf Bounds for the Grundy chromatic number of graphs in terms of domination number}

\author{\large\bf Abbas Khaleghi~~~~~~Manouchehr Zaker\footnote{Corresponding author: mzaker@iasbs.ac.ir}
	\vspace{5mm}\\
	Department of Mathematics,\\
	Institute for Advanced Studies in Basic Sciences,\\
	Zanjan 45137-66731, Iran\\
}

\date{}
\begin{document}
	
\maketitle

\begin{abstract}
\noindent For any graph $G$, the Grundy (or First-Fit) chromatic number of $G$, denoted by $\Gamma(G)$ (also $\chi_{_{\sf FF}}(G)$), is defined as the maximum number of colors used by the First-Fit (greedy) coloring of the vertices of $G$.
Determining the Grundy number is $NP$-complete, and obtaining bounds for $\Gamma(G)$ in terms
of the known graph parameters is an active research topic.
By a star partition of $G$ we mean any partition of $V(G)$ into say $V_1, \ldots, V_k$ such that each $G[V_i]$ contains a vertex adjacent to any other vertex in $V_i$. In this paper using the star partition of graphs we obtain the first upper bounds for the Grundy number in terms of the domination number. We also prove some bounds in terms of the domination number and girth of graphs.
\end{abstract}



\section{Introduction}

\noindent All graphs in this paper are undirected without any loops and multiple edges. Let $G=(V,E)$ be a graph on the vertex set $V=V(G)$ and the edge set $E=E(G)$. By the order and size of $G$ we mean $|V|$ and $|E|$, respectively. Also $\delta(G)$ and $\Delta(G)$ stand for the minimum and the maximum degree of $G$, respectively. For any vertex $v\in V(G)$, denote the neighborhood set of $v$ in $G$ (resp. the closed neighborhood of $v$ in $G$) by $N(v)$ (resp. $N[v]$). For any subset $S$ of vertices in $G$, by $G[S]$ we mean the subgraph of $G$ induced by the elements of $S$. The {\it girth} of $G$ is defined as the smallest length of cycles in $G$. If the graph is acyclic then its girth is infinity. A subset $D\subseteq V(G)$ is said to be a {\it dominating set} of $G$ if any vertex in $V(G)\setminus D$ is adjacent to some vertex of $D$. The {\it domination number} $\gamma(G)$ is the smallest cardinality of any dominating set in $G$. We refer the readers to \cite{BM} for the basic concepts not defined here.

\noindent A {\it Grundy $k$-coloring} of a graph $G$ is a proper $k$-coloring of vertices in $G$ using colors $\{1, 2, \ldots, k\}$ such that for each $i$ and $j$ with $i < j$, any vertex colored $j$ is adjacent to some vertex colored $i$. The {\it First-Fit} or {\it Grundy chromatic number} (or simply Grundy number) of a graph $G$, denoted by $\Gamma(G)$, is the largest integer $k$, such that there exists a Grundy $k$-coloring for $G$. Let $\sigma$ be an ordering of vertices in $G$ and $v_1, \ldots, v_n$ be the order of vertices in $\sigma$. The {\it greedy coloring algorithm} corresponding to $\sigma$, assigns color $1$ to $v_1$ and a minimum available color from positive integers to $v_i$, for any $i\geq 2$. Let $FF(G, \sigma)$ be the number of colors used by the greedy algorithm corresponding to $\sigma$. As proved in \cite{Z2}, $\Gamma(G)=\max FF(G, \sigma)$, where the maximum is taken over all orderings $\sigma$ of $V(G)$. Let $\Gamma(G)=k$. Then any vertex of color $k$ in a Grundy $k$-coloring of $G$ has degree at least $k-1$. We conclude that $\Delta(G) \geq deg(v) \geq k-1$ or $\Gamma(G)\leq \Delta(G)+1$. For each positive integer $k$, a class of graphs denoted by ${\mathcal{A}}_k$ was constructed in \cite{Z2} which satisfies the following property. The Grundy number of any graph $G$ is at least $k$ if and only if $G$ contains an induced subgraph isomorphic to some element of ${\mathcal{A}}_k$. Any element of ${\mathcal{A}}_k$ is called {\it $k$-atom}. For any integer $k\geq 1$, there exists exactly one tree $k$-atom, denoted by $T_{k}$. For $k=1,2$, $T_k$ is isomorphic to the complete graph on one and two vertices, respectively. Assume that $T_k$ is constructed for $k\geq 2$, the tree $T_{k+1}$ is obtained from $T_k$ by attaching one leaf to each vertex of $T_k$ so that $|V(T_{k+1})|=2|V(T_k)|$. Hence, $|V(T_{k})|=2^{k-1}$. For example $T_{3}$ is the path on $4$ vertices. The concept of Grundy number was introduced and studied first time in \cite{CS}. The First-Fit coloring and the Grundy number of graphs were widely studied in graph theory, e.g. \cite{BHLY, GL, Z2}. Also to obtain bounds for the Grundy number of graphs in terms of various graph parameters is the research subject of many researchers, e.g. \cite{C, FGSS, TWHZ, Z3, Z34, Z4}.

\noindent The chromatic number is an $NP$-hard problem even for the graphs with arbitrary large girth \cite{EHK}. Also to determine the Grundy number of graphs of girth at least four is $NP$-hard \cite{HS}. To obtain bounds for the Grundy number in terms of girth is the research subject of many papers e.g. \cite{FGSS, Z3, Z4}. The following result was proved in \cite{Z3}

\begin{thm}(\cite{Z3})
For any graph $G$ whose girth $g$ is an odd integer,
$$\Gamma(G)\leq \frac{g-1}{2}|V(G)|^{\frac{2}{g-1}}.$$\label{zakerbound}
\end{thm}

\noindent As proved in \cite{FGSS}, the bound of Theorem \ref{zakerbound} is sharp for $g=5$. A conjecture made by the second author \cite{Z3} claims that the power $2/(g-1)$ in Theorem \ref{zakerbound} can not be improved. In this paper we improve the bound of Theorem \ref{zakerbound} by replacing the term $|V(G)|$ by $|V(G)|-\gamma(G)$. The Grundy number of graphs containing no cycle of even length $2t$ was investigated in \cite{Z4}, where $t\geq 2$ is fixed and arbitrary integer. Precisely, it was proved in \cite{Z4} that there exists a function $f(t)$ such that $\Gamma(G) \leq f(t) |V(G)|^{1/t}$ whenever $G$ does not contain $C_{2t}$ as ordinary subgraph.

\noindent In a graph $G$, a vertex $u$ is called {\it apex} if $u$ is adjacent to any other vertex of $G$. By a {\it star graph} $S_k$ we mean the complete bipartite graph $K_{1,k}$, i.e. a tree with exactly one internal vertex and $k$ leaves.

\noindent \begin{defin}
\noindent By a star partition of a graph $G$ we mean a partition of $V(G)$ into say $S_1, \ldots, S_k$ such that for each $i$, $G[S_i]$ contains an apex vertex. Denote by $s(G)$ the smallest value $k$ such that $G$ admits a star partition with $k$ subsets.
\end{defin}

\noindent The following proposition is known in the literature (see e.g. \cite{AFGS}) but we prove it here to make this paper self-contained.

\begin{prop}
For any graph $G$ without isolated vertices, $s(G)=\gamma(G)$.\label{sgamma}
\end{prop}

\noindent \begin{proof}
Let $k=\gamma(G)$ and $D=\{u_1, \ldots, u_k\}$ be any minimum dominating set in $G$. For any $i$, $1\leq i \leq k$, $N(u_i)\not= \varnothing$, since $\delta(G)\geq 1$; also $N(u_i)\setminus D \not= \varnothing$, since otherwise $D\setminus \{u_i\}$ is dominating set. Define $S_1=N[u_1]\setminus D$, $S_2=N[u_2]\setminus (S_1\cup D)$ and in general $S_i=N[u_i]\setminus (S_1\cup \ldots \cup S_{i-1} \cup D)$, for each $3\leq i \leq k$. Consider a collection of subsets ${\mathcal{S}}= \{S_i:~S_i\not= \varnothing, 1\leq i \leq k\}$. By the construction, for any $i$ and $j$ with $S_i, S_j \in {\mathcal{S}}$ we have $S_i \cap S_j = \varnothing$. Any vertex of $V(G)$ belongs to one member of ${\mathcal{S}}$ because $D$ is dominating set. We conclude that ${\mathcal{S}}$ is a star partition of $V(G)$. Also $|{\mathcal{S}}|\leq k$. Then $s(G)\leq |{\mathcal{S}}| \leq \gamma(G)$. We prove now $s(G)\geq \gamma(G)$. For this purpose, let $\{S_1, \ldots, S_t\}$ be any star partition of $G$ with $t=s(G)$. Let $u_i$ be an apex vertex in $S_i$ for each $i$. Clearly, $\{u_1, \ldots, u_t\}$ is a dominating set in $G$.
\end{proof}

\noindent Since the domination number of graphs is an $NP$-complete problem \cite{GJ}, then by Proposition \ref{sgamma} to determine the star partition number of graphs is $NP$-complete too. Despite that the domination and star partition numbers are equal quantities, they are conceptually different. Because in one side, we have a subset of vertices satisfying a certain property but on the other side we have a partition of the whole vertex set satisfying some properties. In proving the bounds of this paper in terms of the domination number, the star partitions are very helpful and provide a convenient proof methodology. For this reason, we use the terminology and concept of star partitions. The star partition number has also connections with partitions into bicliques. Let $G$ be a bipartite graph. Any complete bipartite subgraph of $G$ is called a biclique of $G$. Denote by $c(G)$ the minimum number of vertex disjoint bicliques which cover the vertices of $G$. We have $c(G)\leq s(G)$. We will not prove it here, but it can be shown that $s(G)-c(G)$ can be arbitrarily large for bipartite graphs $G$.

\noindent {\bf The paper is organized as follows.} In Sections $2$ we obtain upper bounds for the Grundy number of general and triangle-free graphs. In Section $3$, the bounds are in terms of the domination and girth of graphs, where the girth is odd number. In particular, we prove (Corollary \ref{cor1}) in this section that if $\Delta(G) \leq (g-1)/2$ then $\Gamma(G)\leq \log (n-s)+2$, where $n$, $s$ and $g$ are the order, domination number and the girth of $G$, respectively. In Section $4$, the bounds are in terms of the domination and girth of graphs, where the girth is even number. We conclude that any exact value or lower bound for the domination number of a graph $G$ yields in a new upper bound for $\Gamma(G)$. Fortunately, the literature is full of lower bound results for the domination number, e.g. \cite{DPW, HHJ, HHS}.

\section{Some general bounds}

\noindent The first result presents an upper bound for the Grundy number of general graphs in terms of the star partition number.

\begin{prop}
\noindent For any graph $G$ on $n$ vertices,

\noindent (i) $\Gamma(G)\leq n - \gamma(G)+1$.
\newline (ii) Equality holds in $(i)$ for a graph $G$ if and only if $V(G)$ is partitioned into two subsets say $Q$ and $D$ such that $G[Q]$ is a complete subgraph on $\Gamma(G)-1$ vertices and $G[D]$ is an independent dominating set in $G$ containing $\gamma(G)$ vertices.
\label{gamma}
\end{prop}

\noindent \begin{proof}
Set for simplicity $\Gamma(G)=k$. To prove $(i)$, let $C$ be a Grundy coloring of $G$ using $k$ colors. Let $v$ be a vertex of color $k$ in $C$. Set $H=G[N[v]]$. We have $|V(H)|\geq k$. The subgraph $H$ is a star subgraph of $G$. Consider a star partition of $G$ consisting of $H$ and $n-k$ subgraphs each containing a single vertex of $V(G)\setminus V(H)$. By the definition, $s(G)\leq n-k+1$. Hence, $k\leq n-\gamma(G)+1$.

\noindent To prove $(ii)$, assume first that $\Gamma(G) = n - \gamma(G)+1$. Set $\Gamma(G)=k$. Let $C_1, \ldots, C_k$ be the color classes in a Grundy $k$-coloring of $G$. We have $\gamma(G)\leq |C_1|$ since $C_1$ is a maximal independent set. Using $k = n - \gamma(G)+1$ and $n=|C_1|+ \cdots + |C_k|$ we obtain $|C_t|=1$ for any $t\geq 2$. It follows that the subgraph induced on $Q= C_2 \cup \ldots \cup C_k$ is complete subgraph on $k-1$ vertices. It is easily seen that $|C_1|=\gamma(G)$. Define $D=C_1$. Now, $V(G)= Q\cup D$ is the desired partition.

\noindent Assume now that $V(G)= Q\cup D$ and $Q$ and $D$ have the aforementioned properties. We have $\Gamma(G)=|Q|+1$ and $\gamma(G)=|D|$. Hence, $|V(G)|=|Q|+|D|=\Gamma(G)-1+\gamma(G)$, as desired. It is easily seen that infinitely many graphs can be constructed satisfying $(ii)$.
\end{proof}

\noindent The next result concerns triangle-free graphs. A graph is triangle-free if it does not contain any triangle.

\begin{thm}
Let $G$ be a triangle-free graph on $n$ vertices. Then
$$\Gamma(G) \leq \frac{n-\gamma(G)+4}{2}.$$\label{triangle-free}
\end{thm}

\noindent\begin{proof}
Let $\Gamma(G)=k$ and $C$ be a Grundy coloring of $G$ using $k$ colors. Let $u$ be a vertex of color $k$ in $C$ and $v_1, \ldots, v_{k-1}$ a sequence of its neighbors such that $v_i$ has color $i$ in $C$ for each $i$. Since the coloring $C$ is a Grundy coloring, then $v_{k-1}$ needs $k-2$ neighbors say $u_1, \ldots, u_{k-2}$ such that $u_i$ has color $i$ in $C$ for each $i$. Note that since $G$ is triangle-free then $\{v_1, \ldots, v_{k-1}\}\cap \{u_1, \ldots, u_{k-2}\} =\varnothing$. Define a subgraph $H$ of $G$ induced on $\{u, v_1, \ldots, v_{k-1}, u_1, \ldots, u_{k-2}\}$. Consider the following partition of $V(H)$ into two star subsets $A$ ad $B$, where $A=\{u, v_1, v_2, \ldots, v_{k-2}\}$ and $B=\{v_{k-1}, u_1, u_2, \ldots, u_{k-2}\}$.

\noindent Set $s=\gamma(G)=s(G)$. Let $S'$ be a star partition of $G$ consisting of the stars $A$, $B$ and $n-(2k-2)$ remaining single vertices in $V(G)\setminus V(H)$. Let $|S'|=s'$. We have $s'=2+n-(2k-2)=n-2k+4$. Since $s$ is the star partition number of $G$, then $s\leq n-2k+4$ or $2k\leq n-s+4$. Hence, $$\Gamma(G) \leq \frac{n-\gamma(G)+4}{2}.$$
\end{proof}

\noindent In the following we show that for each $n$, there exists a graph on $n$ vertices for which the equality holds in the bound of Theorem \ref{triangle-free}.

\begin{prop}
For any integer $n\geq 3$, there exists a triangle-free graph $G$ on $n$ vertices such that
$$\Gamma(G)=\lfloor \frac{n+2}{2}\rfloor=\frac{|V(G)|-\gamma(G)+4}{2}.$$\label{=triangle-free}
\end{prop}

\noindent\begin{proof}
First, let $n=2t \geq 4$ be an even integer. Define $G$ as a complete bipartite graph $K_{t,t}$ from which a matching of size $t-1$ is removed. Assume that the partite sets of $G$ are $X$ and $Y$. Assign the colors $t+1, t-1, t-2, \ldots, 2, 1$ to the vertices of $X$ and assign the colors $t, t-1, t-2, \ldots, 2, 1$ to the vertices of $Y$ in such a way that the two vertices having colors $t+1$ and $t$ are adjacent in $G$. Note that $\Gamma(G)$ is $t+1$, since $\Delta(G)=t$ and $\Gamma(G)\leq \Delta(G)+1$ and also the above-mentioned labels provide a Grundy coloring of $G$ using $t+1$ colors. It is easily seen that the star partition number of $G$ is $2$. We have now $\Gamma(G)=t+1=(2t-2+4)/2$.

\noindent Now, let $n=2t+1$, $t \geq 3$. We construct $G'$ as follows. First, let $H$ be a complete bipartite graph $K_{t,t}$ minus a matching of size $t-1$. We color one part say $X$ of this graph with $t+1, t-1,t-2, \ldots, 2, 1$ and other part say $Y$ with $t, t-1, t-2, \ldots, 2, 1$, so that the vertices with colors $t+1$ and $t$ are adjacent. Then, we add one extra vertex say $w$ of color $1$ to the part $X$ and remove the edge between the vertex of color $1$ in $X$ (other than $w$) and the vertex of color $2$ in $Y$. We finally put an edge between $w$ and the vertex of color $2$ in $Y$. Denote the resulting graph by $G'$, which is illustrated in Figure \ref{knn}. We have $\Delta(G')=t$. The above-mentioned labels provide a Grundy coloring of $G$ using $t+1$ colors. Since $\Gamma(G)\leq \Delta(G)+1$. Hence, $\Gamma(G)=t+1$. Also, the star partition number of $G$ is $3$. Then,
$\Gamma(G)=t+1=(2t+1-3+4)/2$.

\begin{figure}[h]
\begin{center}
\begin{tikzpicture}
\filldraw (5,3) circle (2pt) node[xshift=8pt, yshift=8pt, scale=1pt]{t+1};
\filldraw (6,3) circle (2pt) node[xshift=1pt, yshift=8pt, scale=1pt]{t-1};
\filldraw (7,3) circle (2pt) node[xshift=1pt, yshift=8pt, scale=1pt]{t-2};
\filldraw (10,3) circle (2pt) node[xshift=1pt, yshift=8pt, scale=1pt]{2};
\filldraw (11,3) circle (2pt) node[xshift=1pt, yshift=8pt, scale=1pt]{1};
\filldraw[red] (12,3) circle (2pt) node[color=red, xshift=1pt, yshift=8pt, scale=1pt]{1};

\filldraw (5,1) circle (2pt) node[xshift=1pt, yshift=-8pt, scale=1pt]{t};
\filldraw (6,1) circle (2pt) node[xshift=1pt, yshift=-8pt, scale=1pt]{t-1};
\filldraw (7,1) circle (2pt) node[xshift=1pt, yshift=-8pt, scale=1pt]{t-2};
\filldraw (10,1) circle (2pt) node[xshift=1pt, yshift=-8pt, scale=1pt]{2};
\filldraw (11,1) circle (2pt) node[xshift=1pt, yshift=-8pt, scale=1pt]{1};

\filldraw (8,3) circle (0pt);
\filldraw (9,3) circle (0pt);
\filldraw (8,1) circle (0pt);
\filldraw (9,1) circle (0pt);

\draw (5,3)--(5,1);
\draw (5,3)--(6,1);
\draw (5,3)--(7,1);
\draw (5,3)--(10,1);
\draw (5,3)--(11,1);
\draw (6,3)--(5,1);
\draw (6,3)--(7,1);
\draw (6,3)--(10,1);
\draw (6,3)--(11,1);
\draw (7,3)--(5,1);
\draw (7,3)--(6,1);
\draw (7,3)--(10,1);
\draw (7,3)--(11,1);
\draw (10,3)--(5,1);
\draw (10,3)--(6,1);
\draw (10,3)--(7,1);
\draw (10,3)--(11,1);
\draw (11,3)--(5,1);
\draw (11,3)--(6,1);
\draw (11,3)--(7,1);
\draw (11,3)--(10,1);
\draw[red] (12,3)--(10,1);
\draw[dashed] (8,3)--(9,3);
\draw[dashed] (8,1)--(9,1);
\end{tikzpicture}
\caption{The graph $G'$ in the proof of Proposition \ref{=triangle-free}}\label{knn}
\end{center}
\end{figure}
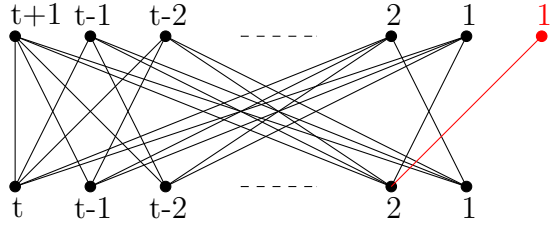

\end{proof}

\noindent In this paragraph we compare the bound of Theorem \ref{triangle-free} with some of the previously known bounds concerning the Grundy number. It was proved in \cite{TWHZ} that if $G$ is triangle-free on $n$ vertices then $\Gamma(G) \leq (n+2)/2$. Note that if $G$ is not isomorphic to a star graph then the bound of Theorem \ref{triangle-free} is better than $(n+2)/2$ bound. In fact, when $G$ is not isomorphic to a star graph then $s(G) \geq 2$ and hence $(n-s(G)+4)/2 \leq (n+2)/2$. The other well-known bound is $\Gamma(G) \leq \Delta(G)+1$. For triangle-free graphs $G$ with $|V(G)|-2\Delta(G)+2 \leq \gamma(G)$, the bound of Theorem \ref{triangle-free} is better than $\Delta(G)+1$ bound. Many triangle-free graphs satisfy this property, as obtained from the lower bounds for the domination number presented in the book \cite{HHS} (Chapter $2$, Theorems 2.22 and 2.28).

\section{Bounds involving the girth of graphs I}

\noindent As mentioned earlier, to obtain bounds for the Grundy number in terms of girth is the research subject of many papers. By Theorem \ref{zakerbound}, if $G$ is any graph on $n$ vertices whose girth $g$ is odd then $\Gamma(G)\leq ((g-1)/2)n^{2/(g-1)}$. A conjecture made by the second author \cite{Z3} claims that the power $2/(g-1)$ in the latter inequality can not be improved. The following results improve the latter bound by replacing $n$ by $n-\gamma(G)$. In this section we consider the case where the girth is odd number. The next section deals with even girth.

\begin{thm}
Let $G$ be a graph on $n$ vertices whose girth $g$ is an odd integer. Then
\begin{center}
$\Gamma(G)\leq {\frac{g-1}{2} (n-\gamma(G))^{\frac{2}{g-1}}}+1.$
\end{center}\label{girth}
\end{thm}

\noindent\begin{proof} \noindent First note that if $n=\gamma(G)$ then the inequality becomes $\Gamma(G)\leq 1$ and this holds since $n=\gamma(G)$ if and only if $G$ is a union of isolated vertices. Assume hereafter that $n-\gamma(G)\geq 1$. Write $s=\gamma(G)=s(G)$. If $g=3$ then the assertion reduces to $\Gamma(G)\leq n-s+1$ which holds by Proposition \ref{gamma} $(i)$. Also if $\Gamma(G) \leq (g+1)/2$ then the inequality holds since $n-\gamma(G)\geq 1$ implies that the right hand side of the inequality is at least $(g+1)/2$. Assume hereafter that $\Gamma(G) \geq \max \{(g+3)/2, 4\}$. We divide the proof into two parts.

\noindent {\bf Case 1: $g \equiv{3} \pmod{4}$}

\noindent Let $\Gamma(G)=k$ and $C$ be a Grundy coloring of $G$ using $k$ colors. Let $v$ be a vertex of color $k$ in $C$. The vertex $v$ has neighbors of colors $1, 2, \ldots, k-1$. These neighbors are mutually non-adjacent, because the girth of $G$ is at least four. Define levels $L_1, L_2, \ldots$ of vertices in $G$ as follow. The $L_1$ consists of the single vertex $v$ which is colored $k$. The level $L_2$ consists of a set of $k-1$ mutually non-adjacent neighbors of $v$ with distinct colors $1, 2, \ldots, k-1$. Any vertex in $L_2$ of color say $j$ needs neighbors of colors $1, 2, \ldots, j-1$. We put these new neighbors in $L_3$. We continue this procedure to obtain the other levels until $L_{(g+1)/2}$. Note that each vertex in $L_{(g+1)/2}$ has distance $(g-1)/2$ from $v$. Because the girth of $G$ is $g$, then the subgraph of $G$ induced on $L_1\cup L_2\cup \ldots \cup L_{(g-1)/2}$ is an induced tree of $G$. Define $H=G[L_1\cup \ldots \cup L_{(g+1)/2}]$. The $L_1$ contains $\binom{k-1}{0}$ vertex. The number of vertices in $L_2$ equals to the number of $1$-element subsets of $\{1, 2, \ldots, k-1\}$, i.e. $\binom{k-1}{1}$. Continuing this argument we obtain that $|L_{(g+1)/2}|$ equals to the number of $(\frac{g-1}{2})$-subsets of $\{1, 2, \ldots,k-1\}$. Hence, $|L_{(g+1)/2}|= \binom{k-1}{\frac{g-1}{2}}$. We conclude that
$$|V(H)|=\sum_{i=0}^{\frac{g-1}{2}}\binom{k-1}{i}.$$

\noindent Consider now a star partition of $V(H)$ into star subgraphs as depicted in Figure \ref{123}. Denote this star partition of $H$ by $S'$. We explain the structure of the star partition $S'$ as follows. All stars in $S'$ have apexes only in the odd-labeled levels. Let $L_i$ be any arbitrary odd-labeled level and let $w$ be an arbitrary vertex in $L_i$. Consider $w$ and all children of $w$ in $L_{i+1}$ (if any) as a star in $S'$. All stars of $S'$ are obtained by this method. In the following we prove that the number of stars in $S'$ equals to $$\sum_{i=0}^{\frac{g-3}{4}}\binom{k-1}{2i}.$$

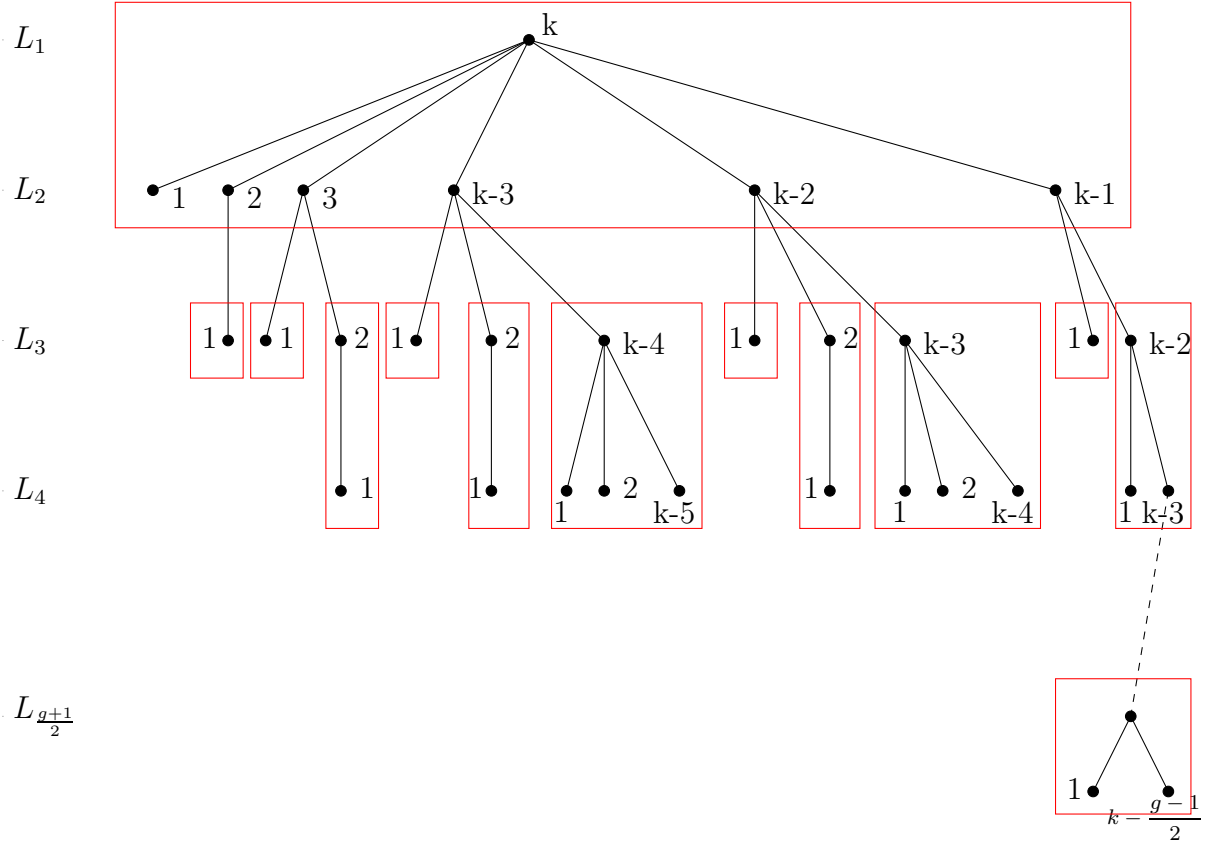
\begin{figure}
\begin{tikzpicture}
\filldraw (8,3) circle (2pt) node[xshift=8pt, yshift=6pt, scale=1pt]{k};
\filldraw (3,1) circle (2pt) node[xshift=10pt, yshift=-3pt]{1};
\filldraw (4,1) circle (2pt) node[xshift=10pt, yshift=-3pt]{2};
\filldraw (5,1) circle (2pt) node[xshift=10pt, yshift=-3pt]{3};
\filldraw (7,1) circle (2pt) node[xshift=15pt, yshift=-1pt]{k-3};
\filldraw (11,1) circle (2pt) node[xshift=15pt, yshift=-1pt]{k-2};
\filldraw (15,1) circle (2pt) node[xshift=15pt, yshift=-1pt]{k-1};
\draw (8,3)--(3,1);
\draw (8,3)--(4,1);
\draw (8,3)--(5,1);
\draw (8,3)--(7,1);
\draw (8,3)--(11,1);
\draw (8,3)--(15,1);

\filldraw (4,-1) circle (2pt) node[xshift=-7pt, yshift=1pt]{1};
\filldraw (4.5,-1) circle (2pt) node[xshift=8pt, yshift=1pt]{1};
\filldraw (5.5,-1) circle (2pt) node[xshift=8pt, yshift=1pt]{2};
\filldraw (6.5,-1) circle (2pt) node[xshift=-7pt, yshift=1pt]{1};
\filldraw (7.5,-1) circle (2pt) node[xshift=8pt, yshift=1pt]{2};
\filldraw (9,-1) circle (2pt) node[xshift=15pt, yshift=-1pt]{k-4};
\filldraw (11,-1) circle (2pt) node[xshift=-7pt, yshift=1pt]{1};
\filldraw (12,-1) circle (2pt) node[xshift=8pt, yshift=1pt]{2};
\filldraw (13,-1) circle (2pt) node[xshift=15pt, yshift=-1pt]{k-3};
\filldraw (15.5,-1) circle (2pt) node[xshift=-7pt, yshift=1pt]{1};
\filldraw (16,-1) circle (2pt) node[xshift=15pt, yshift=-1pt]{k-2};
\draw (4,1)--(4,-1);
\draw (5,1)--(4.5,-1);
\draw (5,1)--(5.5,-1);
\draw (7,1)--(6.5,-1);
\draw (7,1)--(7.5,-1);
\draw (7,1)--(9,-1);
\draw (11,1)--(11,-1);
\draw (11,1)--(12,-1);
\draw (11,1)--(13,-1);
\draw (15,1)--(15.5,-1);
\draw (15,1)--(16,-1);

\filldraw (5.5,-3) circle (2pt) node[xshift=10pt, yshift=1pt]{1};
\filldraw (7.5,-3) circle (2pt) node[xshift=-6pt, yshift=1pt]{1};
\filldraw (8.5,-3) circle (2pt) node[xshift=-2pt, yshift=-8pt]{1};
\filldraw (9,-3) circle (2pt) node[xshift=10pt, yshift=1pt]{2};
\filldraw (10,-3) circle (2pt) node[xshift=-2pt, yshift=-8pt]{k-5};
\filldraw (12,-3) circle (2pt) node[xshift=-7pt, yshift=1pt]{1};
\filldraw (13,-3) circle (2pt) node[xshift=-2pt, yshift=-8pt]{1};
\filldraw (13.5,-3) circle (2pt) node[xshift=10pt, yshift=1pt]{2};
\filldraw (14.5,-3) circle (2pt) node[xshift=-2pt, yshift=-8pt]{k-4};
\filldraw (16,-3) circle (2pt) node[xshift=-2pt, yshift=-8pt]{1};
\filldraw (16.5,-3) circle (2pt) node[xshift=-2pt, yshift=-8pt]{k-3};
\draw (5.5,-1)--(5.5,-3);
\draw (7.5,-1)--(7.5,-3);
\draw (9,-1)--(8.5,-3);
\draw (9,-1)--(9,-3);
\draw (9,-1)--(10,-3);
\draw (12,-1)--(12,-3);
\draw (13,-1)--(13,-3);
\draw (13,-1)--(13.5,-3);
\draw (13,-1)--(14.5,-3);
\draw (16,-1)--(16,-3);
\draw (16,-1)--(16.5,-3);

\filldraw (16,-6) circle (2pt);
\filldraw (15.5,-7) circle (2pt) node[xshift=-7pt, yshift=1pt]{1};
\filldraw (16.5,-7) circle (2pt) node[xshift=-5pt, yshift=-10pt]{\scriptsize ${k-\frac{g-1}{2}}$};
\draw (16,-6)--(15.5,-7);
\draw (16,-6)--(16.5,-7);
\draw[dashed] (16.5,-3)--(16,-6);
\draw[{red}] (2.5,3.5) rectangle (16,.5);
\draw [{red}] (3.5,-0.5) rectangle (4.2,-1.5);
\draw [{red}] (4.3,-0.5) rectangle (5,-1.5);
\draw [{red}] [{red}] (5.3,-0.5) rectangle (6,-3.5);
\draw [{red}] (6.1,-0.5) rectangle (6.8,-1.5);
\draw [{red}] (7.2,-0.5) rectangle (8,-3.5);
\draw [{red}] (8.3,-0.5) rectangle (10.3,-3.5);
\draw [{red}] (10.6,-0.5) rectangle (11.3,-1.5);
\draw [{red}] (11.6,-0.5) rectangle (12.4,-3.5);
\draw [{red}] (12.6,-0.5) rectangle (14.8,-3.5);
\draw [{red}] (15,-0.5) rectangle (15.7,-1.5);
\draw [{red}] (15.8,-0.5) rectangle (16.8,-3.5);
\draw [{red}] (15,-5.5) rectangle (16.8,-7.3);

\filldraw (1,3) circle (0pt) node[right]{$L_{1}$};
\filldraw (1,1) circle (0pt) node[right]{$L_{2}$};
\filldraw (1,-1) circle (0pt) node[right]{$L_{3}$};
\filldraw (1,-3) circle (0pt) node[right]{$L_{4}$};
\filldraw (1,-6) circle (0pt) node[right]{$L_{\frac{g+1}{2}}$};
\end{tikzpicture}
\caption{The subgraph $H$ in Case 1 of the proof of Theorem \ref{girth} with a vertex partition into star subgraphs}\label{123}
\end{figure}

\noindent All stars in $S'$ has apexes only in odd-labeled levels of $H$, as illustrated in Figure \ref{123}. There is one star whose apex is in the first level. The number of stars with apexes in $L_3$ is equal to the number of vertices in $L_3$, i.e. $\binom{k-1}{2}$. The number of stars with apexes in $L_5$ is equal to the number of vertices in $L_5$, i.e. $\binom{k-1}{4}$ and so on. Since $(g-1)/2$ is odd then the number of stars with apexes in $L_{(g-1)/2}$ is  $|L_{(g-1)/2}|=\binom{k-1}{\frac{g-3}{2}}$.

\noindent Let $S''$ be the star partition of $G$ including the stars in $S'$ and the remaining vertices in $V(G)\setminus V(H)$, i.e. each vertex as a single vertex star. Set $|S''|=s''$. Since $s$ is the star partition number of $G$, then
$$s \leq s''=|S'|+|V(G)|-|V(H)| \Longrightarrow |V(H)|-|S'|\leq n-s.$$
\noindent It follows by the previous values for $|V(H)|$ and $|S'|$ that
$$ \binom{k-1}{1}+\binom{k-1}
{3}+\cdots +\binom{k-1}{\frac{g-1}{2}} \leq n-s.$$
Then
$$\binom{k-1}{\frac{g-1}{2}}\leq n-s.$$
\noindent For any two integers $a$ and $b$ with $a \leq b$, we have $(a/b)^{b} \leq \binom{a}{b}$.\label{inequal} Using this inequality we obtain
$$(\frac{k-1}{\frac{g-1}{2}})^{\frac{g-1}{2}}\leq \binom{k-1}{\frac{g-1}{2}}\leq n-s.$$
\noindent It follows that
$$\frac{2}{g-1}(k-1)\leq (n-s)^{\frac{2}{g-1}}.$$
\noindent This completes the proof of Case 1.

\noindent {\bf Case 2: $g \equiv{1} \pmod{4}$}

\noindent The proof is similar to the proof of Case 1. Let $\Gamma(G)=k$ and $C$ be a Grundy coloring of $G$ using $k$ colors. Let the sets $L_1, L_2, \ldots, L_{(g+1)/2}$ and the subset $H$ be as in the proof of Case 1. We consider a star partition $S'$ for $H$ different from the one given in Case 1. Figure \ref{picgirth2} illustrates the star partition $S'$ of $V(H)$. We explain the structure of the star partition $S'$ as follows. The vertex of color $k$ with its only neighbor of color $1$ forms one star in $S'$. All other stars of $S'$ have apexes only in the even-labeled levels. Let $L_i$ be any arbitrary even-labeled level and let $w$ be an arbitrary vertex in $L_i$. Consider $w$ and all children of $w$ in $L_{i+1}$ (if any) as a star in $S'$. All stars of $S'$ are obtained by this method. We prove that the number of stars in $S'$ equals to

\begin{center}
$\sum_{i=1}^{\frac{g-1}{4}}\binom{k-1}{2i-1}.$
\end{center}

\noindent There are $|L_2|=\binom{k-1}{1}$ stars with apexes in $L_{2}$. There are $|L_4|=\binom{k-1}{3}$ stars with apexes in $L_{4}$ and
so on. Since $(g-1)/2$ is even then the number of stars with apexes in $L_{(g-1)/2}$ is $|L_{(g-1)/2}|=\binom{k-1}{\frac{g-3}{2}}$.

\noindent Let $S''$ be the star partition of $G$ including the stars in $S'$ and the remaining vertices in $V(G)\setminus V(H)$, i.e. each vertex as a single vertex star. Set $|S''|=s''$. Since $s$ is the star partition number of $G$, then
$$\sum_{j=0}^{\frac{g-1}{4}}\binom{k-1}{2j} \leq n-s.$$

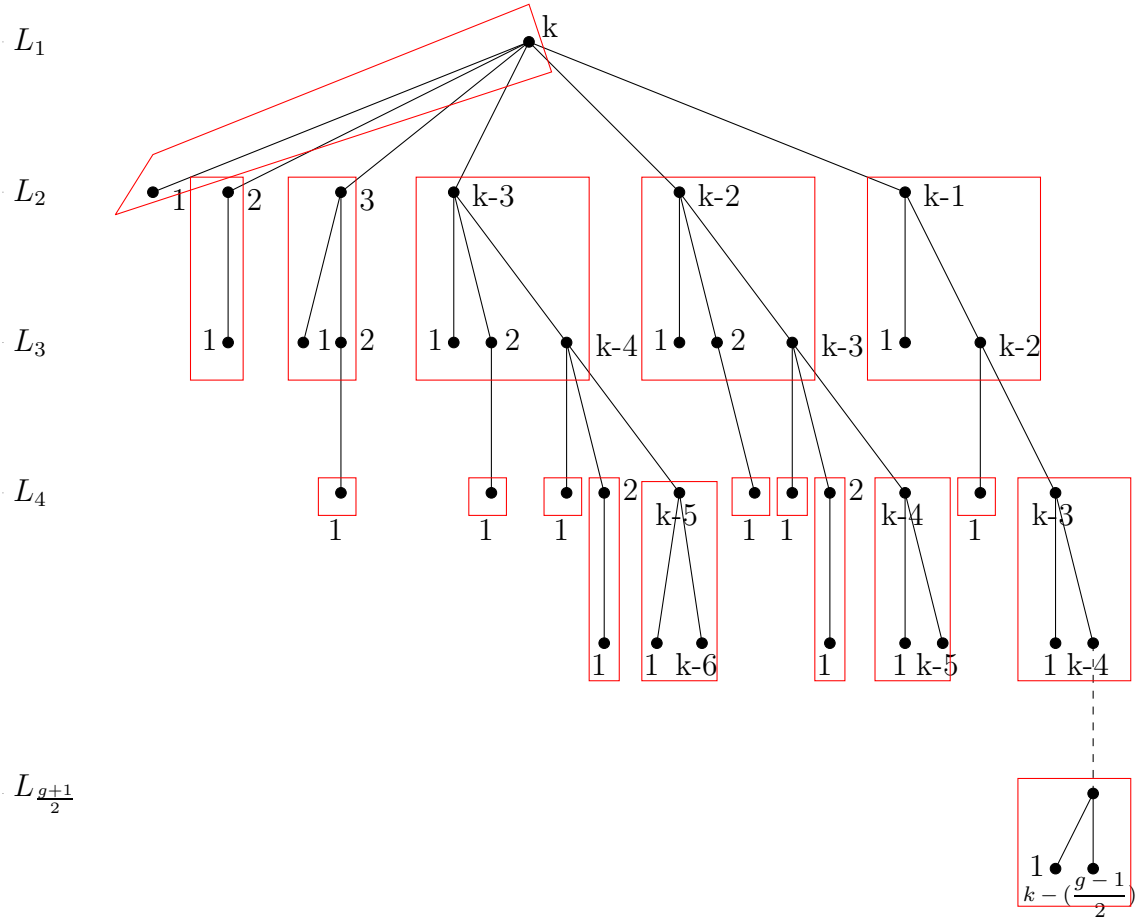
\begin{figure}
\begin{tikzpicture}
\filldraw (8,3) circle (2pt) node[xshift=8pt, yshift=6pt, scale=1pt]{k};

\filldraw (3,1) circle (2pt) node[xshift=10pt, yshift=-3pt]{1};
\filldraw (4,1) circle (2pt) node[xshift=10pt, yshift=-3pt]{2};
\filldraw (5.5,1) circle (2pt) node[xshift=10pt, yshift=-3pt]{3};
\filldraw (7,1) circle (2pt) node[xshift=15pt, yshift=-1pt]{k-3};
\filldraw (10,1) circle (2pt) node[xshift=15pt, yshift=-1pt]{k-2};
\filldraw (13,1) circle (2pt) node[xshift=15pt, yshift=-1pt]{k-1};
\draw (8,3)--(3,1);
\draw (8,3)--(4,1);
\draw (8,3)--(5.5,1);
\draw (8,3)--(7,1);
\draw (8,3)--(10,1);
\draw (8,3)--(13,1);

\filldraw (4,-1) circle (2pt) node[xshift=-7pt, yshift=1pt]{1};
\filldraw (5,-1) circle (2pt) node[xshift=8pt, yshift=1pt]{1};
\filldraw (5.5,-1) circle (2pt) node[xshift=10pt, yshift=1pt]{2};
\filldraw (7,-1) circle (2pt) node[xshift=-7pt, yshift=1pt]{1};
\filldraw (7.5,-1) circle (2pt) node[xshift=8pt, yshift=1pt]{2};
\filldraw (8.5,-1) circle (2pt) node[xshift=19pt, yshift=-1pt]{k-4};
\filldraw (10,-1) circle (2pt) node[xshift=-7pt, yshift=1pt]{1};
\filldraw (10.5,-1) circle (2pt) node[xshift=8pt, yshift=1pt]{2};
\filldraw (11.5,-1) circle (2pt) node[xshift=19pt, yshift=-1pt]{k-3};
\filldraw (13,-1) circle (2pt) node[xshift=-7pt, yshift=1pt]{1};
\filldraw (14,-1) circle (2pt) node[xshift=15pt, yshift=-1pt]{k-2};
\draw (4,1)--(4,-1);
\draw (5.5,1)--(5,-1);
\draw (5.5,1)--(5.5,-1);
\draw (7,1)--(7,-1);
\draw (7,1)--(7.5,-1);
\draw (7,1)--(8.5,-1);
\draw (10,1)--(10,-1);
\draw (10,1)--(10.5,-1);
\draw (10,1)--(11.5,-1);
\draw (13,1)--(13,-1);
\draw (13,1)--(14,-1);

\filldraw (5.5,-3) circle (2pt) node[xshift=-2pt, yshift=-14pt]{1};
\filldraw (7.5,-3) circle (2pt) node[xshift=-2pt, yshift=-14pt]{1};
\filldraw (8.5,-3) circle (2pt) node[xshift=-2pt, yshift=-14pt]{1};
\filldraw (9,-3) circle (2pt) node[xshift=10pt, yshift=1pt]{2};
\filldraw (10,-3) circle (2pt) node[xshift=-1pt, yshift=-8pt]{k-5};
\filldraw (11,-3) circle (2pt) node[xshift=-2pt, yshift=-14pt]{1};
\filldraw (11.5,-3) circle (2pt) node[xshift=-2pt, yshift=-14pt]{1};
\filldraw (12,-3) circle (2pt) node[xshift=10pt, yshift=1pt]{2};
\filldraw (13,-3) circle (2pt) node[xshift=-1pt, yshift=-8pt]{k-4};
\filldraw (14,-3) circle (2pt) node[xshift=-2pt, yshift=-14pt]{1};
\filldraw (15,-3) circle (2pt) node[xshift=-1pt, yshift=-8pt]{k-3};
\draw (5.5,-1)--(5.5,-3);
\draw (7.5,-1)--(7.5,-3);
\draw (8.5,-1)--(8.5,-3);
\draw (8.5,-1)--(9,-3);
\draw (8.5,-1)--(10,-3);
\draw (10.5,-1)--(11,-3);
\draw (11.5,-1)--(11.5,-3);
\draw (11.5,-1)--(12,-3);
\draw (11.5,-1)--(13,-3);
\draw (14,-1)--(14,-3);
\draw (14,-1)--(15,-3);

\filldraw(9,-5) circle (2pt) node[xshift=-2pt, yshift=-8pt]{1};
\filldraw (9.7,-5) circle (2pt) node[xshift=-2pt, yshift=-8pt]{1};
\filldraw (10.3,-5) circle (2pt) node[xshift=-2pt, yshift=-8pt]{k-6};
\filldraw (12,-5) circle (2pt) node[xshift=-2pt, yshift=-8pt]{1};
\filldraw (13,-5) circle (2pt) node[xshift=-2pt, yshift=-8pt]{1};
\filldraw (13.5,-5) circle (2pt) node[xshift=-2pt, yshift=-8pt]{k-5};
\filldraw (15,-5) circle (2pt) node[xshift=-2pt, yshift=-8pt]{1};
\filldraw (15.5,-5) circle (2pt) node[xshift=-2pt, yshift=-8pt]{k-4};

\draw (9,-3)--(9,-5);
\draw (10,-3)--(9.7,-5);
\draw (10,-3)--(10.3,-5);
\draw (12,-3)--(12,-5);
\draw (13,-3)--(13,-5);
\draw (13,-3)--(13.5,-5);
\draw (15,-3)--(15,-5);
\draw (15,-3)--(15.5,-5);

\filldraw (15.5,-7) circle (2pt);
\filldraw (15,-8) circle (2pt) node[xshift=-7pt, yshift=1pt]{1};
\filldraw (15.5,-8) circle (2pt) node[xshift=-5pt, yshift=-10pt]{\scriptsize $k-(\frac{g-1}{2})$};
\draw (15.5,-7)--(15,-8);
\draw (15.5,-7)--(15.5,-8);
\draw [dashed] (15.5,-5)--(15.5,-7);

\draw [{red}] (2.5,0.7)--(8.3,2.6)--(8,3.5)--(3,1.5)--(2.5,0.7);
\draw [{red}] (3.5,1.2) rectangle (4.2,-1.5);
\draw [{red}] (4.8,1.2) rectangle (5.7,-1.5);
\draw [{red}] (6.5,1.2) rectangle (8.8,-1.5);
\draw [{red}] (9.5,1.2) rectangle (11.8,-1.5);
\draw [{red}] (12.5,1.2) rectangle (14.8,-1.5);
\draw [{red}] (5.2,-3.3) rectangle (5.7,-2.8);
\draw [{red}] (7.2,-3.3) rectangle (7.7,-2.8);
\draw [{red}] (8.2,-3.3) rectangle (8.7,-2.8);
\draw [{red}] (8.8,-2.8) rectangle (9.2,-5.5);
\draw [{red}] (9.5,-2.85) rectangle (10.5,-5.5);
\draw [{red}] (10.7,-3.3) rectangle (11.2,-2.8);
\draw [{red}] (11.3,-3.3) rectangle (11.7,-2.8);
\draw [{red}] (11.8,-2.8) rectangle (12.2,-5.5);
\draw [{red}] (12.6,-2.8) rectangle (13.6,-5.5);
\draw [{red}] (13.7,-3.3) rectangle (14.2,-2.8);
\draw [{red}] (14.5,-2.8) rectangle (16,-5.5);
\draw [{red}] (14.5,-6.8) rectangle (16,-8.5);

\filldraw (1,3) circle (0pt) node[right]{$L_{1}$};
\filldraw (1,1) circle (0pt) node[right]{$L_{2}$};
\filldraw (1,-1) circle (0pt) node[right]{$L_{3}$};
\filldraw (1,-3) circle (0pt) node[right]{$L_{4}$};
\filldraw (1,-7) circle (0pt) node[right]{$L_{\frac{g+1}{2}}$};
\end{tikzpicture}
\caption{The subgraph $H$ in Case 2 of the proof of Theorem \ref{girth} with a vertex partition into star subgraphs}\label{picgirth2}
\end{figure}

\noindent We continue the proof exactly the same as the proof of Case 1 but using the star partition $S'$ of $H$ and obtain
$$(\frac{2}{g-1}(k-1))^{\frac{g-1}{2}}\leq n-s.$$
This completes the whole proof.

\end{proof}

\noindent In the following result we used the inequality $\Gamma(G)\leq \Delta(G)+1$.

\begin{cor}
Let $G$ be a graph on $n$ vertices and of odd girth $g$ such that $\Delta(G) \leq (g-1)/2$. Then
\begin{center}
$\Gamma(G)\leq \log (n-\gamma(G))+2.$
\end{center}\label{cor1}
\end{cor}

\noindent\begin{proof}
We prove the corollary for the case $g \equiv{3}\pmod{4}$. The proof for the case $g \equiv{1}\pmod{4}$ is completely similar. Let $g \equiv{3}\pmod{4}$ and $C$ and $L_1, L_2, \ldots, L_{((g-1)/2)+1}$ be as in the proof of Theorem \ref{girth}. Since $\Gamma(G) \leq (g+1)/2$ then the subgraph of $G$ induced on $L_1 \cup L_2 \cup \ldots \cup L_{((g+1)/2)}$ is isomorphic to the tree atom $T_k$, where $k=\Gamma(G)$. We have $|V(T_k)|=2^{k-1}$. Also as demonstrated before $|V(H)|-|S'|\leq n-\gamma(G)$. These imply
$$2^{k-1}-\sum_{i=0}^{\frac{g-3}{4}}\binom{k-1}{2i}\leq n-\gamma(G)$$
and hence
$$2^{k-1}-(1/2)2^{k-1}=2^{k-2}\leq n-\gamma(G).$$
Finally,
$k\leq \log (n-\gamma(G)) +2$.
\end{proof}

\section{Bounds involving the girth of graphs II}

\noindent In this section we obtain bounds for $\Gamma(G)$ involving girth when the girth $g$ is an even integer.
The best known bound for the Grundy number of such graphs is of the form $f(g)|V(G)|^{2/(g-2)}$, for some function $f(g)$ \cite{Z3}. Theorem \ref{girth3} improves this upper bound by replacing $n$ by $n-\gamma(G)$.

\begin{thm}
Let $G$ be a graph on $n$ vertices whose girth $g$ is even. Then
\begin{center}
$\Gamma(G)\leq \frac{g-2}{2}~ (\frac{n-\gamma(G)}{2})^{\frac{2}{g-2}}+2.$
\end{center}\label{girth3}
\end{thm}

\noindent\begin{proof} \noindent If $n-\gamma(G)\leq 1$ then clearly $\Delta(G)\leq 1$ and $G$ has at most one connected component of maximum degree one. The inequality obviously holds for such graphs. Assume hereafter that $n-\gamma(G)\geq 2$. Also if $\Gamma(G) \leq (g+2)/2$ then the inequality holds since by $n-\gamma(G)\geq 2$ the right hand side of the inequality is at least $(g+2)/2$. Assume hereafter that $\Gamma(G) > (g+2)/2$. Write $s=\gamma(G)=s(G)$ throughout the proof. We divide the proof into two parts.

\noindent {\bf Case 1:} $g \equiv{0} \pmod{4}$.

\noindent Let $\Gamma(G)=k$ and $C$ be a Grundy coloring of $G$ using $k$ colors. Let $u$ be a vertex of color $k$ in $C$. The vertex $u$ has a set of neighbors of colors $1, 2, \ldots, k-1$. These neighbors are mutually non-adjacent, because the girth of $G$ is at least four. Let $w$ be a neighbor of $u$ of color $k-1$. Define levels $L_1, L_2, \ldots$ of vertices in $G$ as follow. The $L_1$ consists of two vertices $u$ and $w$. The level $L_2$ consists of a set of $k-2$ mutually non-adjacent neighbors of $u$ with distinct colors $1, 2, \ldots, k-2$ and a set of $k-2$ mutually non-adjacent neighbors of $w$ having distinct colors $1, 2, \ldots, k-2$. Any vertex in $L_2$ of color say $j$ needs neighbors of colors $1, 2, \ldots, j-1$. We put these new neighbors in $L_3$. We continue this procedure to obtain the other levels until $L_{g/2}$. Define $H=G[L_{1}\cup \ldots \cup L_{g/2}]$. Note that $G[L_1 \cup \ldots \cup L_{(g-2)/2}]$ is an induced subtree of $G$ and a possible edge in $H$ is only between two vertices in $L_{g/2}$.
The $L_{1}$ contains $\binom{k-2}{0}+\binom{k-2}{0}$ vertices. The number of vertices in $L_{2}$ equals to twice the number of $1$-element subsets of $\{1, 2, \ldots, k-2\}$, i.e. $\binom{k-2}{1}+\binom{k-2}{1}$. A similar fact holds for the next levels $L_3, \ldots, L_{g/2}$. We conclude that
$$|V(H)|=2\sum_{i=0}^{\frac{g-2}{2}}\binom{k-2}{i}.$$

\noindent Consider now a star partition of $V(H)$ into star subgraphs as depicted in Figure \ref{1234}.

\noindent Denote this star partition of $H$ by $S'$. The structure of $S'$ is explained similar to the previous cases. In the following we prove that the number of stars in $S'$ equals to $$2\sum_{i=0}^{\frac{g-4}{4}}\binom{k-2}{2i}.$$
\noindent There are two stars whose apexes are in the first level. The number of stars with apexes in $L_3$ is equal to the number of vertices in $L_3$, i.e. $\binom{k-2}{2}+\binom{k-2}{2}$ and so on. Since $g/2$ is even then the number of stars with apexes in $L_{(g-2)/2}$ is $|L_{(g-2)/2}|=\binom{k-2}{\frac{g-4}{2}}+\binom{k-2}{\frac{g-4}{2}}$.
\noindent Let $S''$ be the star partition of $G$ including the stars in $S'$ and the remaining vertices in $V(G)\setminus V(H)$, i.e. each vertex as a single vertex star. Set $|S''|=s''$. Since $s$ is the star partition number of $G$, then
$$s \leq s''=|S'|+|V(G)|-|V(H)| \Longrightarrow |V(H)|-|S'|\leq n-s.$$
\noindent By replacing the values of $|V(H)|$ and $|S'|$ in the latter inequality we obtain
$$ \binom{k-2}{1}+\binom{k-2}{3}+\cdots +\binom{k-2}{\frac{g-2}{2}} \leq\frac{1}{2}(n-s).$$
\noindent Then
$$\binom{k-2}{\frac{g-2}{2}}\leq \frac{1}{2}(n-s).$$

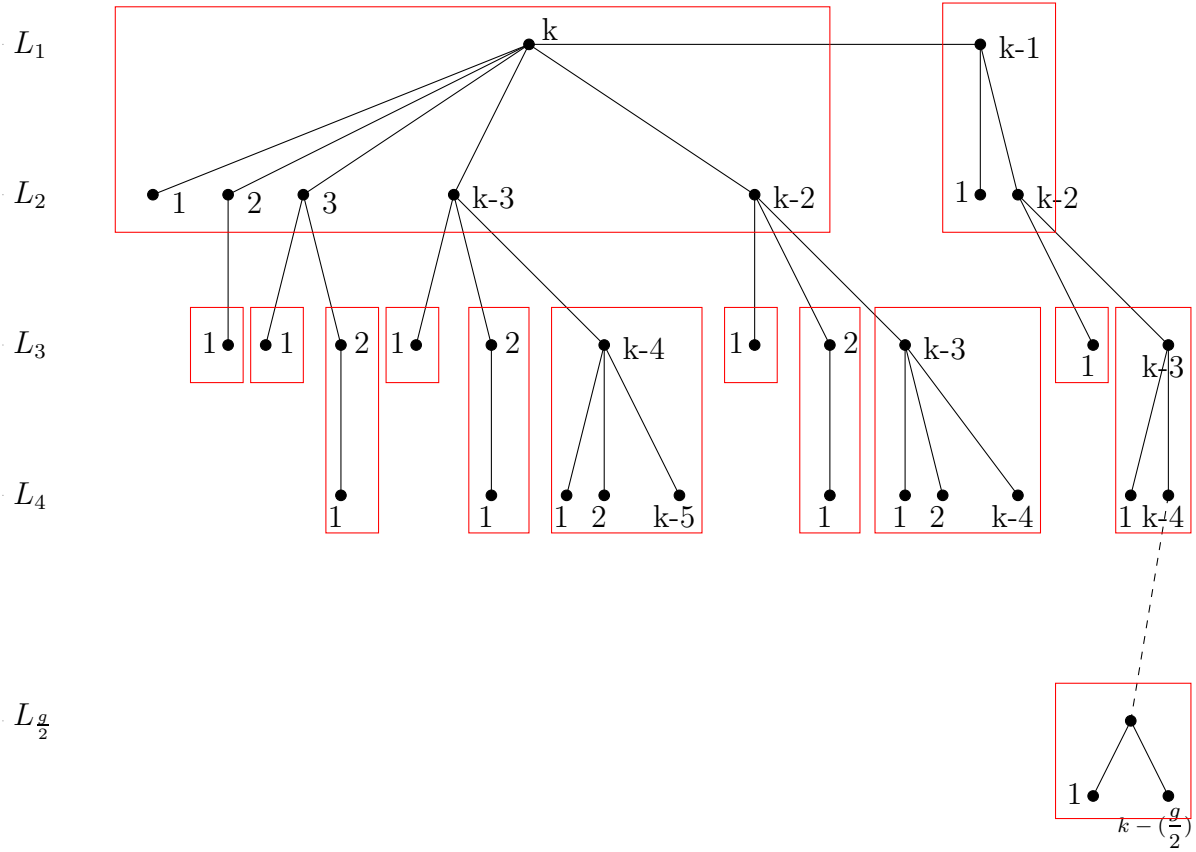
\begin{figure}
\begin{tikzpicture}
\filldraw (8,3) circle (2pt) node[xshift=8pt, yshift=6pt, scale=1pt]{k};
\filldraw (3,1) circle (2pt) node[xshift=10pt, yshift=-3pt]{1};
\filldraw (4,1) circle (2pt) node[xshift=10pt, yshift=-3pt]{2};
\filldraw (5,1) circle (2pt) node[xshift=10pt, yshift=-3pt]{3};
\filldraw (7,1) circle (2pt) node[xshift=15pt, yshift=-1pt]{k-3};
\filldraw (11,1) circle (2pt) node[xshift=15pt, yshift=-1pt]{k-2};
\filldraw (14,3) circle (2pt) node[xshift=15pt, yshift=-1pt]{k-1};
\draw (8,3)--(3,1);
\draw (8,3)--(4,1);
\draw (8,3)--(5,1);
\draw (8,3)--(7,1);
\draw (8,3)--(11,1);
\draw (8,3)--(14,3);

\filldraw (4,-1) circle (2pt) node[xshift=-7pt, yshift=1pt]{1};
\filldraw (4.5,-1) circle (2pt) node[xshift=8pt, yshift=1pt]{1};
\filldraw (5.5,-1) circle (2pt) node[xshift=8pt, yshift=1pt]{2};
\filldraw (6.5,-1) circle (2pt) node[xshift=-7pt, yshift=1pt]{1};
\filldraw (7.5,-1) circle (2pt) node[xshift=8pt, yshift=1pt]{2};
\filldraw (9,-1) circle (2pt) node[xshift=15pt, yshift=-1pt]{k-4};
\filldraw (11,-1) circle (2pt) node[xshift=-7pt, yshift=1pt]{1};
\filldraw (12,-1) circle (2pt) node[xshift=8pt, yshift=1pt]{2};
\filldraw (13,-1) circle (2pt) node[xshift=15pt, yshift=-1pt]{k-3};
\filldraw (14,1) circle (2pt) node[xshift=-7pt, yshift=1pt]{1};
\filldraw (14.5,1) circle (2pt) node[xshift=15pt, yshift=-1pt]{k-2};
\draw (4,1)--(4,-1);
\draw (5,1)--(4.5,-1);
\draw (5,1)--(5.5,-1);
\draw (7,1)--(6.5,-1);
\draw (7,1)--(7.5,-1);
\draw (7,1)--(9,-1);
\draw (11,1)--(11,-1);
\draw (11,1)--(12,-1);
\draw (11,1)--(13,-1);
\draw (14,3)--(14,1);
\draw (14,3)--(14.5,1);

\filldraw (5.5,-3) circle (2pt) node[xshift=-2pt, yshift=-8pt]{1};
\filldraw (7.5,-3) circle (2pt) node[xshift=-2pt, yshift=-8pt]{1};
\filldraw (8.5,-3) circle (2pt) node[xshift=-2pt, yshift=-8pt]{1};
\filldraw (9,-3) circle (2pt) node[xshift=-2pt, yshift=-8pt]{2};
\filldraw (10,-3) circle (2pt) node[xshift=-2pt, yshift=-8pt]{k-5};
\filldraw (12,-3) circle (2pt) node[xshift=-2pt, yshift=-8pt]{1};
\filldraw (13,-3) circle (2pt) node[xshift=-2pt, yshift=-8pt]{1};
\filldraw (13.5,-3) circle (2pt) node[xshift=-2pt, yshift=-8pt]{2};
\filldraw (14.5,-3) circle (2pt) node[xshift=-2pt, yshift=-8pt]{k-4};
\filldraw (16,-3) circle (2pt) node[xshift=-2pt, yshift=-8pt]{1};
\filldraw (16.5,-3) circle (2pt) node[xshift=-2pt, yshift=-8pt]{k-4};
\filldraw (15.5,-1) circle (2pt) node[xshift=-2pt, yshift=-8pt]{1};
\filldraw (16.5,-1) circle (2pt) node[xshift=-2pt, yshift=-8pt]{k-3};
\draw (5.5,-1)--(5.5,-3);
\draw (7.5,-1)--(7.5,-3);
\draw (9,-1)--(8.5,-3);
\draw (9,-1)--(9,-3);
\draw (9,-1)--(10,-3);
\draw (12,-1)--(12,-3);
\draw (13,-1)--(13,-3);
\draw (13,-1)--(13.5,-3);
\draw (13,-1)--(14.5,-3);
\draw (14.5,1)--(15.5,-1);
\draw (14.5,1)--(16.5,-1);
\draw (16.5,-1)--(16,-3);
\draw (16.5,-1)--(16.5,-3);

\filldraw (16,-6) circle (2pt);
\filldraw (15.5,-7) circle (2pt) node[xshift=-7pt, yshift=1pt]{1};
\filldraw (16.5,-7) circle (2pt) node[xshift=-5pt, yshift=-12pt]{\scriptsize $k-(\frac{g}{2})$};
\draw (16,-6)--(15.5,-7);
\draw (16,-6)--(16.5,-7);
\draw[dashed] (16.5,-3)--(16,-6);
\draw [{red}] (2.5,3.5) rectangle (12,.5);
\draw [{red}] (13.5,3.55) rectangle (15,0.5);
\draw [{red}] (3.5,-0.5) rectangle (4.2,-1.5);
\draw [{red}] (4.3,-0.5) rectangle (5,-1.5);
\draw [{red}] (5.3,-0.5) rectangle (6,-3.5);
\draw [{red}] (6.1,-0.5) rectangle (6.8,-1.5);
\draw [{red}] (7.2,-0.5) rectangle (8,-3.5);
\draw [{red}] (8.3,-0.5) rectangle (10.3,-3.5);
\draw [{red}] (10.6,-0.5) rectangle (11.3,-1.5);
\draw [{red}] (11.6,-0.5) rectangle (12.4,-3.5);
\draw [{red}] (12.6,-0.5) rectangle (14.8,-3.5);
\draw [{red}] (15,-0.5) rectangle (15.7,-1.5);
\draw [{red}] (15.8,-0.5) rectangle (16.8,-3.5);
\draw [{red}] (15,-5.5) rectangle (16.8,-7.3);

\filldraw (1,3) circle (0pt) node[right]{$L_{1}$};
\filldraw (1,1) circle (0pt) node[right]{$L_{2}$};
\filldraw (1,-1) circle (0pt) node[right]{$L_{3}$};
\filldraw (1,-3) circle (0pt) node[right]{$L_{4}$};
\filldraw (1,-6) circle (0pt) node[right]{$L_{\frac{g}{2}}$};
\end{tikzpicture}
\caption{The subgraph $H$ in Case 1 of the proof of Theorem \ref{girth3} with a vertex partition into star subgraphs}\label{1234}
\end{figure}

\noindent Applying $(a/b)^{b} \leq \binom{a}{b}$ for the latter inequality we obtain the following which completes the proof.
$$(\frac{2}{g-2}(k-2))^{\frac{g-2}{2}}\leq \frac{1}{2}(n-s).$$

\noindent This completes the proof of Case 1.

\noindent {\bf Case 2: $g \equiv{2} \pmod{4}$}.

\noindent The proof is similar to the proof of Case 1. Let $\Gamma(G)=k$ and $C$ be a Grundy coloring of $G$ using $k$
colors. Let the sets $L_1, L_2, \ldots, L_{g/2}$ and the subset $H$ be as in the proof of Case 1. We consider a star partition $S'$ for $H$ different from the one given in Case 1. Figure \ref{12345} illustrates the star partition $S'$ of $V(H)$. We prove that the number of stars in $S'$ equals to
\begin{center}
$2\sum_{i=1}^{\frac{g-2}{4}}\binom{k-2}{2i-1}.$
\end{center}

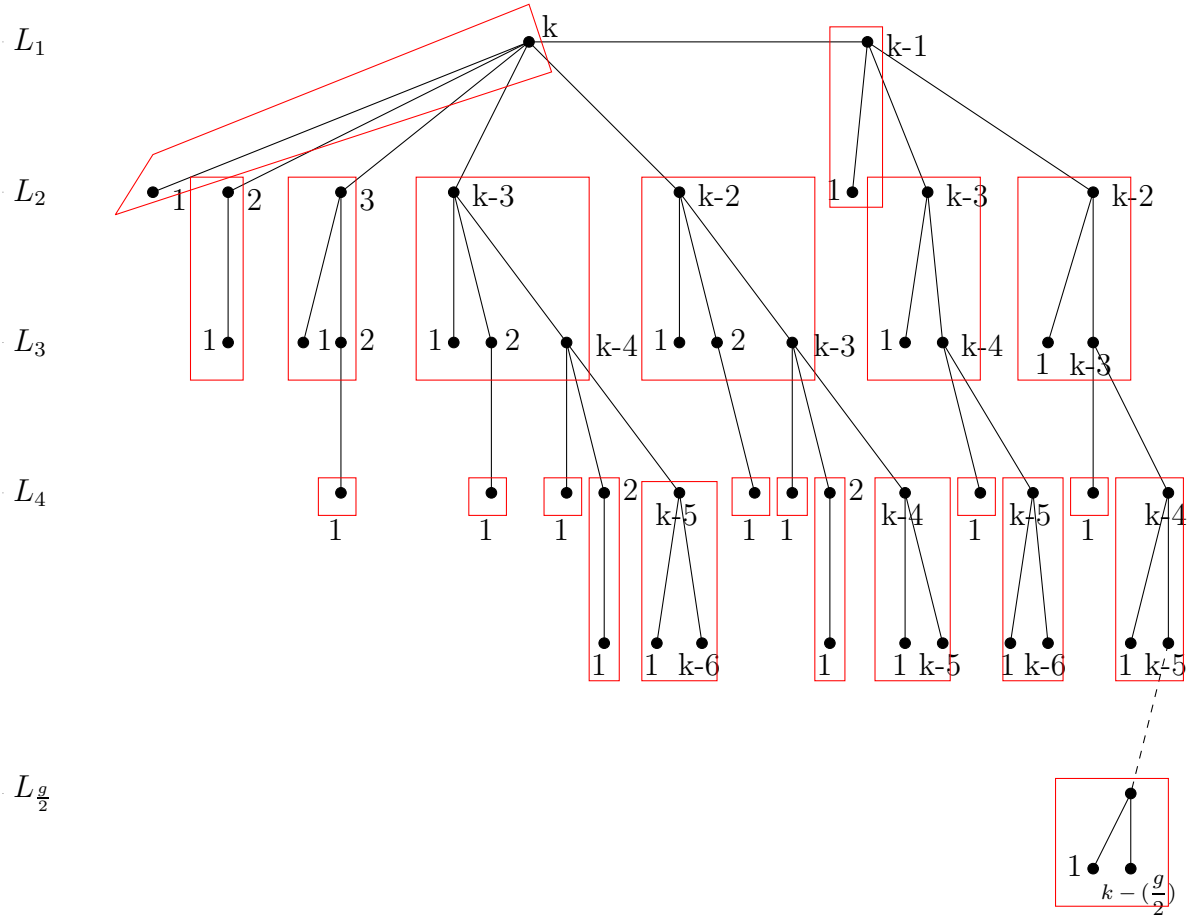
\begin{figure}
\begin{tikzpicture}
\filldraw (8,3) circle (2pt) node[xshift=8pt, yshift=6pt, scale=1pt]{k};

\filldraw (3,1) circle (2pt) node[xshift=10pt, yshift=-3pt]{1};
\filldraw (4,1) circle (2pt) node[xshift=10pt, yshift=-3pt]{2};
\filldraw (5.5,1) circle (2pt) node[xshift=10pt, yshift=-3pt]{3};
\filldraw (7,1) circle (2pt) node[xshift=15pt, yshift=-1pt]{k-3};
\filldraw (10,1) circle (2pt) node[xshift=15pt, yshift=-1pt]{k-2};
\filldraw (13.3,1) circle (2pt) node[xshift=15pt, yshift=-1pt]{k-3};
\filldraw (12.5,3) circle (2pt) node[xshift=15pt, yshift=-1pt]{k-1};
\draw (8,3)--(3,1);
\draw (8,3)--(4,1);
\draw (8,3)--(5.5,1);
\draw (8,3)--(7,1);
\draw (8,3)--(10,1);
\draw (8,3)--(12.5,3);
\draw (12.5,3)--(13.3,1);
\filldraw (4,-1) circle (2pt) node[xshift=-7pt, yshift=1pt]{1};
\filldraw (5,-1) circle (2pt) node[xshift=8pt, yshift=1pt]{1};
\filldraw (5.5,-1) circle (2pt) node[xshift=10pt, yshift=1pt]{2};
\filldraw (7,-1) circle (2pt) node[xshift=-7pt, yshift=1pt]{1};
\filldraw (7.5,-1) circle (2pt) node[xshift=8pt, yshift=1pt]{2};
\filldraw (8.5,-1) circle (2pt) node[xshift=19pt, yshift=-1pt]{k-4};
\filldraw (10,-1) circle (2pt) node[xshift=-7pt, yshift=1pt]{1};
\filldraw (10.5,-1) circle (2pt) node[xshift=8pt, yshift=1pt]{2};
\filldraw (11.5,-1) circle (2pt) node[xshift=16pt, yshift=-1pt]{k-3};
\filldraw (12.3,1) circle (2pt) node[xshift=-7pt, yshift=1pt]{1};
\filldraw (15.5,1) circle (2pt) node[xshift=15pt, yshift=-1pt]{k-2};

\filldraw (13,-1) circle (2pt) node[xshift=-7pt, yshift=1pt]{1};
\filldraw (13.5,-1) circle (2pt) node[xshift=15pt, yshift=-1pt]{k-4};

\draw (4,1)--(4,-1);
\draw (5.5,1)--(5,-1);
\draw (5.5,1)--(5.5,-1);
\draw (7,1)--(7,-1);
\draw (7,1)--(7.5,-1);
\draw (7,1)--(8.5,-1);
\draw (10,1)--(10,-1);
\draw (10,1)--(10.5,-1);
\draw (10,1)--(11.5,-1);
\draw (12.5,3)--(12.3,1);
\draw (12.5,3)--(15.5,1);

\draw (13.3,1)--(13,-1);
\draw(13.3,1)--(13.5,-1);

\filldraw (5.5,-3) circle (2pt) node[xshift=-2pt, yshift=-14pt]{1};
\filldraw (7.5,-3) circle (2pt) node[xshift=-2pt, yshift=-14pt]{1};
\filldraw (8.5,-3) circle (2pt) node[xshift=-2pt, yshift=-14pt]{1};
\filldraw (9,-3) circle (2pt) node[xshift=10pt, yshift=1pt]{2};
\filldraw (10,-3) circle (2pt) node[xshift=-1pt, yshift=-8pt]{k-5};
\filldraw (11,-3) circle (2pt) node[xshift=-2pt, yshift=-14pt]{1};
\filldraw (11.5,-3) circle (2pt) node[xshift=-2pt, yshift=-14pt]{1};
\filldraw (12,-3) circle (2pt) node[xshift=10pt, yshift=1pt]{2};
\filldraw (13,-3) circle (2pt) node[xshift=-1pt, yshift=-8pt]{k-4};
\filldraw (14.9,-1) circle (2pt) node[xshift=-2pt, yshift=-8pt]{1};
\filldraw (15.5,-1) circle (2pt) node[xshift=-1pt, yshift=-8pt]{k-3};

\filldraw (14,-3) circle (2pt) node[xshift=-2pt, yshift=-14pt]{1};
\filldraw (14.7,-3) circle (2pt) node[xshift=-1pt, yshift=-8pt]{k-5};

\draw (5.5,-1)--(5.5,-3);
\draw (7.5,-1)--(7.5,-3);
\draw (8.5,-1)--(8.5,-3);
\draw (8.5,-1)--(9,-3);
\draw (8.5,-1)--(10,-3);
\draw (10.5,-1)--(11,-3);
\draw (11.5,-1)--(11.5,-3);
\draw (11.5,-1)--(12,-3);
\draw (11.5,-1)--(13,-3);
\draw (15.5,1)--(14.9,-1);
\draw (15.5,1)--(15.5,-1);

\draw (13.5,-1)--(14,-3);
\draw (13.5,-1)--(14.7,-3);

\filldraw(9,-5) circle (2pt) node[xshift=-2pt, yshift=-8pt]{1};
\filldraw (9.7,-5) circle (2pt) node[xshift=-2pt, yshift=-8pt]{1};
\filldraw (10.3,-5) circle (2pt) node[xshift=-1pt, yshift=-8pt]{k-6};
\filldraw (12,-5) circle (2pt) node[xshift=-2pt, yshift=-8pt]{1};
\filldraw (13,-5) circle (2pt) node[xshift=-2pt, yshift=-8pt]{1};
\filldraw (13.5,-5) circle (2pt) node[xshift=-1pt, yshift=-8pt]{k-5};
\filldraw (15.5,-3) circle (2pt) node[xshift=-2pt, yshift=-14pt]{1};
\filldraw (16.5,-3) circle (2pt) node[xshift=-1pt, yshift=-8pt]{k-4};
\filldraw (16,-5) circle (2pt) node[xshift=-2pt, yshift=-8pt]{1};
\filldraw (16.5,-5) circle (2pt) node[xshift=-1pt, yshift=-8pt]{k-5};

\filldraw (14.4,-5) circle (2pt) node[xshift=-1pt, yshift=-8pt]{1};
\filldraw (14.9,-5) circle (2pt) node[xshift=-1pt, yshift=-8pt]{k-6};

\draw (9,-3)--(9,-5);
\draw (10,-3)--(9.7,-5);
\draw (10,-3)--(10.3,-5);
\draw (12,-3)--(12,-5);
\draw (13,-3)--(13,-5);
\draw (13,-3)--(13.5,-5);
\draw (15.5,-1)--(15.5,-3);
\draw (15.5,-1)--(16.5,-3);
\draw (16.5,-3)--(16,-5);
\draw (16.5,-3)--(16.5,-5);
\draw (14.7,-3)--(14.4,-5);
\draw (14.7,-3)--(14.9,-5);

\filldraw (16,-7) circle (2pt);
\filldraw (15.5,-8) circle (2pt) node[xshift=-7pt, yshift=1pt]{1};
\filldraw (16,-8) circle (2pt) node[xshift=3pt, yshift=-10pt]{\scriptsize $k-(\frac{g}{2})$};
\draw (16,-7)--(15.5,-8);
\draw (16,-7)--(16,-8);
\draw [dashed] (16.5,-5)--(16,-7);

\draw [{red}] (2.5,0.7)--(8.3,2.6)--(8,3.5)--(3,1.5)--(2.5,0.7);
\draw [{red}] (12,3.2) rectangle (12.7,0.8);
\draw [{red}] (3.5,1.2) rectangle (4.2,-1.5);
\draw [{red}] (4.8,1.2) rectangle (5.7,-1.5);
\draw [{red}] (6.5,1.2) rectangle (8.8,-1.5);
\draw [{red}] (9.5,1.2) rectangle (11.8,-1.5);
\draw [{red}] (12.5,1.2) rectangle (14,-1.5);
\draw [{red}] (14.5,1.2) rectangle (16,-1.5);
\draw [{red}] (5.2,-3.3) rectangle (5.7,-2.8);
\draw [{red}] (7.2,-3.3) rectangle (7.7,-2.8);
\draw [{red}] (8.2,-3.3) rectangle (8.7,-2.8);
\draw [{red}] (8.8,-2.8) rectangle (9.2,-5.5);
\draw [{red}] (9.5,-2.85) rectangle (10.5,-5.5);
\draw [{red}] (10.7,-3.3) rectangle (11.2,-2.8);
\draw [{red}] (11.3,-3.3) rectangle (11.7,-2.8);
\draw [{red}] (11.8,-2.8) rectangle (12.2,-5.5);
\draw [{red}] (12.6,-2.8) rectangle (13.6,-5.5);
\draw [{red}] (13.7,-3.3) rectangle (14.2,-2.8);
\draw [{red}] (14.3,-2.8) rectangle (15.1,-5.5);
\draw [{red}] (15.2,-3.3) rectangle (15.7,-2.8);
\draw [{red}] (15.8,-2.8) rectangle (16.7,-5.5);
\draw [{red}] (15,-6.8) rectangle (16.5,-8.5);

\filldraw (1,3) circle (0pt) node[right]{$L_{1}$};
\filldraw (1,1) circle (0pt) node[right]{$L_{2}$};
\filldraw (1,-1) circle (0pt) node[right]{$L_{3}$};
\filldraw (1,-3) circle (0pt) node[right]{$L_{4}$};
\filldraw (1,-7) circle (0pt) node[right]{$L_{\frac{g}{2}}$};
\end{tikzpicture}
\caption{The subgraph $H$ in Case 2 of the proof of Theorem \ref{girth3} with a vertex partition into star subgraphs}\label{12345}
\end{figure}

\noindent There are $|L_2|=2\binom{k-2}{1}$ stars with apexes in $L_{2}$. There are $|L_4|=2\binom{k-2}{3}$ stars with apexes in $L_{4}$ and
so on. Since $(g-2)/2$ is even then the number of stars with apexes in $L_{(g-2)/2}$ is $|L_{(g-2)/2}|=2\binom{k-2}{\frac{g-4}{2}}$.
With a similar argument as in the proof of Theorem \ref{girth3}, we have
$$|V(H)|=2\sum_{i=0}^{\frac{g-2}{2}}\binom{k-2}{i}.$$
\noindent Let $S''$ be the star partition of $G$ including the stars in $S'$ and the remaining vertices in $V(G)\setminus V(H)$, i.e. each vertex as a single vertex star. Set $|S''|=s''$. Since $s$ is the star partition number of $G$, then
$$\sum_{j=0}^{\frac{g-2}{4}}\binom{k-2}{2j} \leq \frac{1}{2}(n-s).$$
We continue the proof exactly the same as the proof of Theorem \ref{girth3} but using the star partition $S'$ of $H$ and obtain.
$$(\frac{2}{g-2}(k-2))^{\frac{g-2}{2}}\leq \frac{1}{2}(n-s).$$
This completes the proof.
\end{proof}



\end{document}